\documentclass[a4paper,11pt]{amsart}
\usepackage{amssymb}

\numberwithin{equation}{section}
\newtheorem{theorem}{Theorem}[section]

\newtheorem{lemma}[theorem]{Lemma}

\theoremstyle{remark}

\author[J.~Benameur]{Jamel Benameur}
\address{Department of Mathematics, College of Science, King Saud University\\
Riyadh 11451, Kingdom of Saudi Arabia}
\email{\sl jbenameur@ksu.edu.sa}
\author[M.~Benhamed]{Moez Benhamed}
\address{University of Tunis El Manar, Faculty of Sciences of Tunis, LR03ES04, 2092 Tunis, Tunisia}
\email{\sl moez.benhamed@hotmail.com}

\title[Global existence of the two-dimensional
QGE with sub-critical dissipation]{Global existence of the two-dimensional
QGE with sub-critical dissipation}
\date{\today}
\begin{document}
\begin{abstract}
In this paper, we study the sub-critical dissipative quasi-geostrophic equations $({\bf S}_\alpha)$. We prove that there exists a unique local-in-time solution for any large initial
data $\theta_0$ in the space ${\bf{\mathcal X}}^{1-2\alpha}(\mathbb R^2)$ defined by (\ref{ec}). Moreover we show that $({\bf S}_\alpha)$ has a global solution in time if the norms of the initial data in ${\bf{\mathcal X}}^{1-2\alpha}(\mathbb R^2)$ are bounded by $1/4$. Also, we prove a blow-up criterion of the non global solution of $({\bf S}_\alpha)$.
\end{abstract}
\maketitle
\tableofcontents
\section{Introduction}
In this article, we study the initial value-problem for the two-dimensional
quasi-geostrophic equation with sub-critical dissipation $({\bf S}_\alpha)\;, \; 1/2<\alpha\leq 1$,
$$\left\{\begin{array}{l}
  \displaystyle\partial_t
\theta+(-\Delta)^{\alpha}\theta +u.\nabla \theta =0 \\
\theta(0)=\theta_0
\end{array}\right.
\leqno({\bf S}_\alpha) $$
where $1/2<\alpha\leq 1$ is a real number. The variable $\theta$ represents potential temperature, $u=(\partial_2(-\Delta)^{-1/2}\theta,\; -\partial_1(-\Delta)^{-1/2}\theta)$ is the fluid velocity. In the following, we are interested in the case when $1/2<\alpha\leq 1$ as the case $\alpha=1/2$ is very much similar to the $2D$ Navier-Stokes which is known to be a wellposed problem, and therefore the
case $\alpha>1/2$ is even easier to deal with.\\
The mathematical study of the non-dissipative case has first been proposed by Constantin, Majda
and Tabak in \cite{CMT} where it is shown to be an analogue to the $3D$ Euler equations. The dissipative case
has then been studied by Constantin and Wu in \cite{CW11} when $\alpha > 1/2$ and global existence in Sobolev
spaces is studied by Constantin.\\

In this paper, we study $({\bf S}_\alpha)$ in scaling invariant spaces. Solvability of evolution equations in scaling invariant spaces is
well-developed in the context of the Navier-Stokes equations. For example, if we restrict the function spaces to the energy
spaces, the optimal result is due to Fujita and Kato in \cite{HFTK}. Later, Chemin \cite{JYC} proved similar results in the framework of Besov
spaces $\dot{B}_{p, q}^{\frac{q}{p}-1}$. Let us find the scaling invariant critical spaces for $({\bf S}_\alpha)$. The equation is invariant under the following scaling:
$$\theta_\lambda(t, x) = \lambda^{2\alpha-1}\theta(\lambda^{2\alpha}t, \lambda x),\quad \mathrm{with\; initial\; data}\;\; \theta_\lambda(0, x) = \theta_\lambda^0(x) = \lambda^{2\alpha-1}\theta^0(\lambda x).$$
So $L^{\frac{2}{2\alpha-1}}, \dot{H}^{2-2\alpha}\;\mathrm{and}\; \dot{B}_{p, q}^{1+\frac{2}{p}-2\alpha}$ are critical spaces. The global well-posedness for small initial data in critical Besov spaces with
$p < \infty$ was obtained in \cite{HB, WU}. The global well-posedness for large data in critical spaces was obtained by several authors;
in Lebesgue space $L^{\frac{2}{2\alpha-1}}$ by Carrillo and Ferreira \cite{JCLF} for $\alpha> \frac{1}{2}$
, in energy space $H^1$ by Dong and Du \cite{HDDD}, and in Besov spaces
by Abidi and Hmidi \cite{HATH}, Hmidi and Keraani \cite{THSK} and Wang and Zhang \cite{WZ} with $p =\infty$.\\

In this paper, we will solve the system $({\bf S}_\alpha)$ in the critical space by means of a contraction argument. Thus, we can obtain a unique local solution to the
system $({\bf S}_\alpha)$ for any initial data in the critical space and prove the corresponding solution
will be global if the initial data is sufficient small. Before giving our main result let us first precise the notion of
critical space.

For $\sigma \in \mathbb R$, we define the functional space
\begin{equation}\label{ec} {\bf{\mathcal X}}^{\sigma}(\mathbb R^2)=\Big\{f\in \mathcal{D'}(\mathbb R^2);\;(\xi\mapsto|\xi|^{\sigma}\widehat{f}(\xi))\in L^1(\mathbb R^2) \Big\}\end{equation}
which  is equipped with the norm
$$\|f\|_{{\bf{\mathcal X}}^{\sigma}}=\int_{\mathbb R^2}|\xi|^{\sigma}|\widehat{f}(\xi)| \;d\xi.$$

Now we give our first result overcoming the above mentioned smallness assumption. Our global existence result reads as follows.
\begin{theorem}\label{thjm1} Let $\theta^0\in{\bf{\mathcal X}}^{1-2\alpha}(\mathbb R^2)$. There is a time $T>0$ and unique solution
$\theta\in\mathcal C([0,T],{\bf{\mathcal X}}^{1-2\alpha}(\mathbb R^2))$ of $({\bf S}_\alpha)$, moreover $\theta\in L^1([0,T],{\bf{\mathcal X}}^{1}(\mathbb R^2))$. If $\|\theta^0\|_{{\bf{\mathcal X}}^{1-2\alpha}}<1/4$, there are global existence and
\begin{equation}\label{thjm1eq}
\|\theta\|_{{\bf{\mathcal X}}^{1-2\alpha}}+(1-4\|\theta^0\|_{{\bf{\mathcal X}}^{1-2\alpha}})\int_0^t\|\theta\|_{{\bf{\mathcal X}}^{1}}\leq \|\theta^0\|_{{\bf{\mathcal X}}^{1-2\alpha}},\;\;\forall t\geq0.
\end{equation}
\end{theorem}
The proof of the above theorem is based on  Fixed Point Theorem. We will establish a local existence result
and will be able to get global existence that is essentially based on energy estimate of the solution in the space
${\bf{\mathcal X}}^{1-2\alpha}(\mathbb R^2)$.\\
 We state now our second main result.
\begin{theorem}\label{thjm2}
Let $\theta\in\mathcal C([0,T^*),{\bf{\mathcal X}}^{1-2\alpha}(\mathbb R^2))$ be the maximal solution of $({\bf S}_\alpha)$ given by Theorem \ref{thjm1}. Then
\begin{equation}\label{thjm2eq}
T^*<\infty\Longrightarrow\int_0^{T^*}\|\theta\|_{{\bf{\mathcal X}}^{1}}=\infty.
\end{equation}
\end{theorem}
The rest of this paper is divided into three sections. Section \ref{S2} is further divided into tow
subsections. Section \ref{S2.1} recalls the notations, Section \ref{S2.2} preliminaries results. Section \ref{S3} is further divided into three
subsections. Section \ref{S3.1} deals with the existence of solution with any large initial data
in the spaces in ${\bf{\mathcal X}}^{1-2\alpha}(\mathbb R^2)$ , section \ref{S3.2} deals with the uniqueness of solution, and section \ref{S3.3} s devoted to  global existence. The section \ref{S4} proves a blow-up criterion of the non global solution given by Theorem \ref{thjm1}.
\section{Notations and Preliminaries results}\label{S2}
\subsection{Notations} \label{S2.1} In this short paragraph, we give some notations\\\\
$\bullet$ For $f$, we denote  $u_f$ the following
$$u_f=(\partial_2(-\Delta)^{-1/2}f,\; -\partial_1(-\Delta)^{-1/2}f)$$
 $\bullet$ The Fourier transform ${\mathcal F(f)}$ of a tempered distribution $f$ on $\mathbb R^2$ is defined as $${\mathcal F(f)}(\xi)=\int_{\mathbb R^2}e^{-ix\xi}f(x)\;dx.$$
       $\bullet$ The inverse Fourier formula is $${\mathcal F^{-1}(f)}(\xi)= \frac{1}{(2\pi)^2}\int_{\mathbb R^2}e^{ix\xi}f(x)\;dx.$$
       $\bullet$ For any Banach space $(B, \|.\|)$, any real number $1\leq p\leq \infty$ and any time $T>0$, we will denote by $L^p_T(B)$ the space of all measurable functions $t\in [0, T]\longmapsto f(t)\in B$ such that $(t\longmapsto \|f(t)\|)\in L^p([0, T]).$\\\\ $\bullet$ The fractional Laplacian operator $(-\Delta)^{\alpha}$ for a real number $\alpha$ is defined through the Fourier transform, namely $$\widehat{(-\Delta)^{\alpha}f}(\xi)=|\xi|^{2\alpha}\widehat{f}(\xi).$$
\subsection{Preliminaries results}\label{S2.2} The main result of this section is the following lemmas that will play a
crucial role in the proof of our main theorem.
\begin{lemma}\label{Lm2.1} Let $\alpha \in ]\frac{1}{2}, 1]$. We have the following inequalities
 \begin{equation}\label{eq34}\|f\|_{{\bf{\mathcal X}}^{2-2\alpha}}\leq \|f\|_{{\bf{\mathcal X}}^{1-2\alpha}}^{1-\frac{1}{2\alpha}}\|f\|_{{\bf{\mathcal X}}^{1}}^{\frac{1}{2\alpha}}\end{equation}
and \begin{eqnarray}\|f\|_{{\bf{\mathcal X}}^{0}}\leq \|f\|_{{\bf{\mathcal X}}^{1-2\alpha}}^{\frac{1}{2\alpha}}\|f\|_{{\bf{\mathcal X}}^{1}}^{1-\frac{1}{2\alpha}}.\label{eq35}\end{eqnarray}
\end{lemma}
{\bf Proof.}\\
$\bullet$ \textbf{Proof of} (\ref{eq34}). \\ By choosing $(p,q)=(2\alpha,\displaystyle\frac{2\alpha}{2\alpha-1})$ and applying the H\"{o}lder's Inequality, we deduce that \begin{eqnarray*}
\|f\|_{{\bf{\mathcal X}}^{2-2\sigma}}&=&\int_{\xi}|\xi|^{2-2\sigma}|\widehat{f}(\xi)| \;d\xi\\
&\leq&\int_{\xi}|\xi|^{\frac{1}{2\alpha}}|\widehat{f}(\xi)|^{\frac{1}{2\alpha}}|\xi|^{\frac{2\alpha-1}{2\alpha}}
|\widehat{f}(\xi)|^{\frac{2\alpha-1}{2\alpha}}|\xi|^{1-2\alpha} \;d\xi\\
&\leq&\Big(\int_{\xi}|\xi||\widehat{f}(\xi)| \;d\xi\Big)^{1/2\alpha}\Big(\int_{\xi}|\xi|^{1-2\sigma}|\widehat{f}(\xi)| \;d\xi\Big)^{1-\frac{1}{2\alpha}}\\
\end{eqnarray*} Hence,
\begin{eqnarray*}
\|f\|_{{\bf{\mathcal X}}^{2-2\sigma}}&\leq&\|f\|_{{\bf{\mathcal X}}^{1}}^{1/2\alpha}\|f\|_{{\bf{\mathcal X}}^{1-2\sigma}}^{1-\frac{1}{2\alpha}}.
\end{eqnarray*}\endproof
\noindent $\bullet$ \textbf{Proof of} (\ref{eq35}).\\
To estimate $\|f\|_{{\bf{\mathcal X}}^{0}}$, we choose $\displaystyle p=2\alpha \;\;\mathrm{and}\;\; q=\frac{2\alpha}{2\alpha-1}$ and apply the H\"{o}lder's Inequalities. We obtain,
\begin{eqnarray*}
\|f\|_{{\bf{\mathcal X}}^{0}}&=&\int_{\xi}|\xi|^{0}|\widehat{f}(\xi)| \;d\xi\\
&\leq&\int_{\xi}|\xi|^{\frac{2\alpha-1}{2\alpha}}|\widehat{f}(\xi)|^{\frac{2\alpha-1}{2\alpha}}|\xi|^{\frac{1-2\alpha}{2\alpha}}
|\widehat{f}(\xi)|^{\frac{1}{2\alpha}}\;d\xi\\
&\leq&\Big(\int_{\xi}|\xi||\widehat{f}(\xi)| \;d\xi\Big)^{1/2\alpha}\Big(\int_{\xi}|\xi|^{1-2\alpha}|\widehat{f}(\xi)| \;d\xi\Big)^{1-\frac{1}{2\alpha}}\\
\end{eqnarray*} Hence,
\begin{eqnarray*}
\|f\|_{{\bf{\mathcal X}}^{0}}&\leq&\|f\|_{{\bf{\mathcal X}}^{1}}^{1-\frac{1}{2\alpha}}\|f\|_{{\bf{\mathcal X}}^{1-2\alpha}}^{\frac{1}{2\alpha}}.
\end{eqnarray*}
The following lemma, which is a direct consequence of the preceding one, will be
useful in the proof of Theorem \ref{thjm1}.
\begin{lemma}\label{Lm2.2}
Let $f,g\in {\bf{\mathcal X}}^{1-2\alpha}(\mathbb R^2)\cap{\bf{\mathcal X}}^{1}(\mathbb R^2)$. We have the following inequalities
\begin{eqnarray}
 \|fg\|_{{\bf{\mathcal X}}^{2-2\alpha}}&\leq& 2\|f\|_{{\bf{\mathcal X}}^{2-2\alpha}}\|g\|_{{\bf{\mathcal X}}^{0}}+2\|f\|_{{\bf{\mathcal X}}^{0}}\|g\|_{{\bf{\mathcal X}}^{2-2\alpha}}.\label{eq36} \end{eqnarray}
\begin{eqnarray}\quad\quad\quad\|fg\|_{{\bf{\mathcal X}}^{2-2\alpha}}&\leq& 2\|f\|_{{\bf{\mathcal X}}^{1-2\alpha}}^{1-\frac{1}{2\alpha}}\|f\|_{{\bf{\mathcal X}}^{1}}^{\frac{1}{2\alpha}}
\|g\|_{{\bf{\mathcal X}}^{1-2\alpha}}^{\frac{1}{2\alpha}}\|g\|_{{\bf{\mathcal X}}^{1}}^{1-\frac{1}{2\alpha}}\nonumber\\
&+&
2\|f\|_{{\bf{\mathcal X}}^{1-2\alpha}}^{\frac{1}{2\alpha}}\|f\|_{{\bf{\mathcal X}}^{1}}^{1-\frac{1}{2\alpha}}
\|g\|_{{\bf{\mathcal X}}^{1-2\alpha}}^{1-\frac{1}{2\alpha}}\|g\|_{{\bf{\mathcal X}}^{1}}^{\frac{1}{2\alpha}}.\label{eq37}
 \end{eqnarray}
 Particularly
 \begin{eqnarray}\|f. f\|_{{\bf{\mathcal X}}^{2-2\alpha}}\leq 4\|f\|_{{\bf{\mathcal X}}^{1-2\alpha}}\|f\|_{{\bf{\mathcal X}}^{1}}.\label{eq40}
 \end{eqnarray}
 \end{lemma}
{\bf Proof.}\\
$\bullet$ {\bf {Proof of (\ref{eq36}).}}\\
\begin{eqnarray*} \|fg\|_{{\bf{\mathcal X}}^{2-2\alpha}}&=&\int_{\xi}|\xi|^{2-2\alpha}|\widehat{fg}(\xi)| \;d\xi\\
&=& \int_{\xi}|\xi|^{2-2\alpha}|\widehat{f}*\widehat{g}(\xi)| \;d\xi \\
&=& \int_{\xi}|\xi|^{2-2\alpha}|\int_{\eta}\widehat{f}(\eta)\widehat{g}(\xi-\eta)| \;d\xi.
\end{eqnarray*} Then
\begin{eqnarray*}
\|fg\|_{{\bf{\mathcal X}}^{2-2\alpha}}&\leq&\int_{\xi}\int_{\eta}|\xi|^{2-2\alpha}|\widehat{f}(\eta)|
|\widehat{g}(\xi-\eta)| \;d\xi\\
&\leq&\int_{\xi}|\xi|^{2-2\alpha}\int_{\eta}|\widehat{f}(\eta)||\widehat{g}(\xi-\eta)| \;d\xi\\
&\leq&I_1+I_2.
\end{eqnarray*}
Where,
\begin{eqnarray*}
&I_1:= \displaystyle\int_{\xi}|\xi|^{2-2\alpha}\int_{|\eta|<|\xi-\eta|}|\widehat{f}(\eta)||\widehat{g}(\xi-\eta)| \;d\xi\\
&I_2:= \displaystyle\int_{\xi}|\xi|^{2-2\alpha}\int_{|\eta|>|\xi-\eta|}|\widehat{f}(\eta)||\widehat{g}(\xi-\eta)| \;d\xi.
\end{eqnarray*}On the one hand, if $|\eta|<|\xi-\eta|$,
we have $$|\xi|\leq |\xi-\eta|+|\eta|\leq 2\max(|\xi-\eta|, |\eta|)=2|\xi-\eta|.$$ Then  \begin{eqnarray*}
|\xi|^{2-2\alpha}&\leq& 2^{2-2\alpha}|\xi-\eta|^{2-2\alpha}\\
&\leq&2|\xi-\eta|^{2-2\alpha}, \;\;\;\mathrm{for\; all}\;\alpha\in (1/2,1].
\end{eqnarray*}Therefore
\begin{eqnarray*}
I_1&=& \int_{\xi}|\xi|^{2-2\alpha}\int_{|\eta|<|\xi-\eta|}|\widehat{f}(\eta)||\widehat{g}(\xi-\eta)| \;d\xi\\
&\leq& 2\int_{\xi}\int_{\eta}|\widehat{f}(\eta)||\xi-\eta|^{2-2\alpha}|\widehat{g}(\xi-\eta)| \;d\xi\\
&\leq&2\||\widehat{f}|*|\xi|^{2-2\alpha}|\widehat{g}|\|_{L^1}.
\end{eqnarray*} By Young inequality, we get
\begin{eqnarray}I_1\leq 2\|f\|_{{\bf{\mathcal X}}^{0}}\|g\|_{{\bf{\mathcal X}}^{2-2\alpha}}.\label{eq38}\end{eqnarray}On the other hand, if $|\eta|>|\xi-\eta|$,
we have $$|\xi|\leq |\xi-\eta|+|\eta|\leq 2\max(|\xi-\eta|, |\eta|)=2|\eta|.$$ Then  \begin{eqnarray*}
|\xi|^{2-2\alpha}&\leq& 2^{2-2\alpha}|\eta|^{2-2\alpha}\\
&\leq&2|\eta|^{2-2\alpha}, \;\;\;\mathrm{for\; all}\;\alpha\in (1/2,1].\end{eqnarray*}Hence
\begin{eqnarray*}
I_2&=& \int_{\xi}|\xi|^{2-2\alpha}\int_{|\eta|>|\xi-\eta|}|\widehat{f}(\eta)||\widehat{g}(\xi-\eta)| \;d\xi\\
&\leq& 2\int_{\xi}\int_{\eta}|\eta|^{2-2\alpha}|\widehat{f}(\eta)||\widehat{g}(\xi-\eta)| \;d\xi\\
&\leq& 2||\xi|^{2-2\alpha}|\widehat{f}|*|\widehat{g}|\|_{L^1}.
\end{eqnarray*} Young inequality gives
\begin{eqnarray}I_2\leq 2\|f\|_{{\bf{\mathcal X}}^{2-2\alpha}}\|g\|_{{\bf{\mathcal X}}^{0}}.\label{eq39}\end{eqnarray}
Combining (\ref{eq38}) and (\ref{eq39}), we get (\ref{eq36}).  \\

$\bullet$ {\bf {Proof of (\ref{eq37}).}}\\
By including the inequalities (\ref{eq34}) and (\ref{eq35}) in (\ref{eq36}), we obtain
$$\|fg\|_{{\bf{\mathcal X}}^{2-2\alpha}}\leq 2\|f\|_{{\bf{\mathcal X}}^{1-2\alpha}}^{1-\frac{1}{2\alpha}}\|f\|_{{\bf{\mathcal X}}^{1}}^{\frac{1}{2\alpha}}
\|g\|_{{\bf{\mathcal X}}^{1-2\alpha}}^{\frac{1}{2\alpha}}\|g\|_{{\bf{\mathcal X}}^{1}}^{1-\frac{1}{2\alpha}}
+
2\|f\|_{{\bf{\mathcal X}}^{1-2\alpha}}^{\frac{1}{2\alpha}}\|f\|_{{\bf{\mathcal X}}^{1}}^{1-\frac{1}{2\alpha}}
\|g\|_{{\bf{\mathcal X}}^{1-2\alpha}}^{1-\frac{1}{2\alpha}}\|g\|_{{\bf{\mathcal X}}^{1}}^{\frac{1}{2\alpha}}.$$ \endproof
\noindent The proof of Theorem \ref{thjm1} requires the following lemma.
\begin{lemma} \label{Lm2.4} Let $\theta\in L^\infty_T({\bf{\mathcal X}}^{1-2\alpha}(\mathbb R^2))\cap L^1_T({\bf{\mathcal X}}^{1}(\mathbb R^2))$. Then
$$\|\int_0^t e^{-(t-z)|D|^{2\alpha}}{\rm div}\,(\theta u_\theta)dz\|_{{\bf{\mathcal X}}^{1-2\alpha}}\leq C_\alpha
\|\theta \|_{L^\infty_T({\bf{\mathcal X}}^{1-2\alpha})}\|\theta\|_{L^1_T({\bf{\mathcal X}}^1)}.$$ \end{lemma}
{\bf Proof.}
 We have
\begin{eqnarray*}
 \|\int_0^t e^{-(t-z)|D|^{2\alpha}}{\rm div}\,(\theta u_\theta)dz\|_{{\bf{\mathcal X}}^{1-2\alpha}}&\leq&\int_0^t \|e^{-(t-z)|D|^{2\alpha}}{\rm div}\,(\theta u_\theta)\|_{{\bf{\mathcal X}}^{1-2\alpha}}dz\\
  &\leq&\int_0^t \|e^{-(t-z)|D|^{2\alpha}}(\theta u_\theta)\|_{{\bf{\mathcal X}}^{2-2\alpha}}dz\\
   &\leq&\int_0^t \|\theta u_\theta\|_{{\bf{\mathcal X}}^{2-2\alpha}}dz.
   \end{eqnarray*} Using Lemma \ref{Lm2.2}, we obtain,
   \begin{eqnarray*}
    \|\int_0^t e^{-(t-z)|D|^{2\alpha}}{\rm div}\,(\theta u_\theta)dz\|_{{\bf{\mathcal X}}^{1-2\alpha}}&\leq&2\int_0^t \|\theta\|_{{\bf{\mathcal X}}^{2-2\alpha}}\|u_\theta\|_{{\bf{\mathcal X}}^0}+\|\theta\|_{{\bf{\mathcal X}}^0}\|u_\theta\|_{{\bf{\mathcal X}}^{2-2\alpha}}dz.
  \end{eqnarray*}
Using Lemma \ref{Lm2.4} and the fact $\|u_\theta\|_{{\bf{\mathcal X}}^{1}}\leq \|\theta\|_{{\bf{\mathcal X}}^1}$ and $\|u_\theta\|_{{\bf{\mathcal X}}^{1-2\alpha}}\leq \|\theta\|_{{\bf{\mathcal X}}^{1-2\alpha}}$, we get
\begin{eqnarray*}
  \|\int_0^t e^{-(t-z)|D|^{2\alpha}}{\rm div}\,(\theta u_\theta)dz\|_{{\bf{\mathcal X}}^{1-2\alpha}}&\leq&4\int_0^t \|\theta\|_{{\bf{\mathcal X}}^{1-2\alpha}}\|\theta\|_{{\bf{\mathcal X}}^1}dz\\
 &\leq&4\|\theta\|_{L^\infty_T({\bf{\mathcal X}}^{1-2\alpha})}\|\theta\|_{L^1_T({\bf{\mathcal X^1}})}.
  \end{eqnarray*}
Also the following lemma will be useful in the sequel.
\begin{lemma}\label{Lm2.5} Let $\theta\in L^\infty_T({\bf{\mathcal X}}^{1-2\alpha}(\mathbb R^2))\cap L^1_T({\bf{\mathcal X}}^{1}(\mathbb R^2))$. Then
$$\int_0^T\|\int_0^t e^{-(t-z)|D|^{2\alpha}}{\rm div}\,(\theta u_\theta)dz\|_{{\bf{\mathcal X}}^{1}}dt\leq 4
\|\theta \|_{L^\infty_T({\bf{\mathcal X}}^{1-2\alpha})}\|\theta\|_{L^1_T({\bf{\mathcal X}}^1)}.$$ \end{lemma}
{\bf Proof.} We have
\begin{eqnarray*}
 \int_0^T\|\int_0^t e^{-(t-z)|D|^{2\alpha}}{\rm div}\,(\theta u_\theta)dz\|_{{\bf{\mathcal X}}^{1}}dt&\leq&\int_0^T\int_0^t \|e^{-(t-z)|D|^{2\alpha}}{\rm div}\,(\theta u_\theta)(z)\|_{{\bf{\mathcal X}}^{1}}dzdt\\
 &\leq&\int_0^T\int_0^t \int_\xi e^{-(t-z)|\xi|^{2\alpha}}|\xi|^2|\mathcal F(\theta u_\theta)(z,\xi)|d\xi dz dt\\
 &\leq&\int_\xi|\xi|^2\Big(\int_0^T\int_0^t  e^{-(t-z)|\xi|^{2\alpha}}|\mathcal F(\theta u_\theta)(z,\xi)|\;dz\; dt\Big)d\xi.
 \end{eqnarray*}
By integrating the function $e^{-(t-z)|\xi|^{2\alpha}}$ twice with respect to $z\in [0,t]$ and with respect to $t\in[0,T]$, we get
 \begin{eqnarray*}
 \int_0^T\int_0^t e^{-(t-z)|\xi|^{2\alpha}}|\mathcal F(\theta u_\theta)(z,\xi)| dz dt=\int_0^T|\mathcal F(\theta u)(z,\xi)|\Big(\int_z^T e^{-(t-z)|\xi|^{2\alpha}}\;dt\Big)\;dz.
 \end{eqnarray*} Then
 \begin{eqnarray*}
 \int_0^T\|\int_0^t e^{-(t-z)|D|^{2\alpha}}{\rm div}\,(\theta u_\theta)dz\|_{{\bf{\mathcal X}}^{1}}dt&\leq&\int_\xi|\xi|^2\Big(\int_0^T\Big[\int_z^T  e^{-(t-z)|\xi|^{2\alpha}}dt\Big]|\mathcal F(\theta u_\theta)(z,\xi)| dz\Big)d\xi\\
  &\leq&\int_\xi|\xi|^2\Big(\int_0^T\Big[\frac{1-e^{-(T-z)|\xi|^{2\alpha}}}{|\xi|^{2\alpha}}\Big]|\mathcal F(\theta u_\theta)(z,\xi)| dz\Big)d\xi\\
  &\leq&\int_\xi|\xi|^{2-2\alpha}\int_0^T|\mathcal F(\theta u_\theta)(z,\xi)| dz\Big)d\xi\\
  &\leq&\int_0^T\|(\theta u_\theta)(z)\|_{{\bf{\mathcal X}}^{2-2\alpha}} dz.
  \end{eqnarray*}
Using Lemma \ref{Lm2.2} and the fact $\|u_\theta\|_{{\bf{\mathcal X}}^{1}}\leq \|\theta\|_{{\bf{\mathcal X}}^1}$ and $\|u_\theta\|_{{\bf{\mathcal X}}^{1-2\alpha}}\leq \|\theta\|_{{\bf{\mathcal X}}^{1-2\alpha}}$, we get
$$\begin{array}{ccc}
 \displaystyle\int_0^T\|\int_0^t e^{-(t-z)|D|^{2\alpha}}{\rm div}\,(\theta u_\theta)dz\|_{{\bf{\mathcal X}}^{1}}\;dt&\leq&\displaystyle4\int_0^T\|\theta (z)\|_{{\bf{\mathcal X}}^{2-2\alpha}}\|\theta (z)\|_{{\bf{\mathcal X}}^1} dz\\
 &\leq&\displaystyle4\|\theta\|_{L^\infty_T({\bf{\mathcal X}}^{1-2\alpha})}\|\theta\|_{L^1_T({\bf{\mathcal X}}^1)}.
  \end{array}$$
\section{Proof of Theorem \ref{thjm1}}\label{S3}
\subsection{Existence} \label{S3.1} The idea of the proof is to write the initial condition as sum higher and lower frequencies. For small frequencies, we give a regular solution of the associated linear system to $ ({\bf S}_\alpha) $ and for the higher frequencies we consider  a partial differential equation very similar to $ ({\bf S}_\alpha) $  with small initial data in $ {\bf {\mathcal X}} ^ {1-2 \alpha} $ for which we can solve it by the fixed point theorem.\\
\noindent$\bullet$ Let $r$ be a real number such that $0<r<1/20$.\\
\noindent$\bullet$ Let $N\in\mathbb N$, such that
$$\int_{|\xi|>N} |\xi|^{1-2\alpha}|\widehat{\theta}^0(\xi)|d\xi<\frac{r}{5}.$$
Put $a^0$ and $b^0$ defined by
$$\begin{array}{ccc}
a^0&=&\mathcal F^{-1}({\bf 1}_{|\xi|<N}\widehat{\theta^0}(\xi))\\
b^0&=&\mathcal F^{-1}({\bf 1}_{|\xi|>N}\widehat{\theta^0}(\xi)).\\
\end{array}$$
Clearly \begin{equation}\label{eqv01}\|b^0\|_{{\bf {\mathcal X}}^{1-2\alpha}}<\frac{r}{5}.\end{equation}
And put
$$a=e^{-t|D|^{2\alpha}}a^0,$$
$a$ is the unique solution of the heat equation
$$\left\{\begin{array}{ccc}
\partial_t a+(-\Delta)^\alpha a&=&0\\
a(0)&=&a^0.
\end{array}\right.$$
We have
$$\|a\|_{{\bf{\mathcal X}}^{1-2\alpha}}\leq\|\theta^0\|_{{\bf {\mathcal X}}^{1-2\alpha}},\;\;\forall t\geq0,$$
and
$$\begin{array}{ccc}
\|a\|_{L^1([0,T],{\bf{\mathcal X}}^{1})}&=&\displaystyle\int_0^T\int_\xi e^{-t|\xi|^{2\alpha}}|\xi||\widehat{\theta^0}(\xi)|d\xi dt\\
&=&\displaystyle\int_\xi\Big(\int_0^T e^{-t|\xi|^{2\alpha}}dt\Big)|\xi||\widehat{\theta^0}(\xi)|d\xi\\
&=&\displaystyle\int_\xi(1- e^{-T|\xi|^{2\alpha}})|\xi|^{1-2\alpha}|\widehat{\theta^0}(\xi)|d\xi.
\end{array}$$
Using Dominated Convergence Theorem, we get
\begin{equation}\label{eqv0}\lim_{T\rightarrow 0^+}\|a\|_{L^1([0,T],{\bf{\mathcal X}}^{1})}=0.\end{equation}
Let $\varepsilon>0$ such that
\begin{equation}\label{eqv1}
4\varepsilon\|\theta^0\|_{{\bf{\mathcal X}}^{1-2\alpha})}<\frac{r}{5}
\end{equation}
and
\begin{equation}\label{eqv2}
2\Big(\|\theta^0\|_{{\bf{\mathcal X}}^{1-2\alpha}}^{1-\frac{1}{2\alpha}}\varepsilon^{\frac{1}{2\alpha}}+
\|\theta^0\|_{{\bf{\mathcal X}}^{1-2\alpha})}^{\frac{1}{2\alpha}}
\varepsilon^{1-\frac{1}{2\alpha}}\Big)<\frac{1}{5}.
\end{equation}
By equation (\ref{eqv0}), there is a time $T=T(\varepsilon)>0$ such that
\begin{equation}\|a\|_{L^1([0,T],{\bf{\mathcal X}}^{1})}<\varepsilon.\end{equation}
Put $b=\theta-a$, clearly $b$ is a solution of the following system
$$\left\{\begin{array}{cc}
&\partial_t b+(-\Delta)^\alpha b+(u_a+u_b)\nabla(a+b)=0\\
&b(0)=b^0.
\end{array}\right.$$
The integral form of $b$ is as follow
$$b=e^{-t|D|^{2\alpha}}b^0-\int_0^te^{-(t-z)|D|^{2\alpha}}(u_a+u_b)\nabla(a+b)dz.$$
To prove the existence of $b$, put the following operator
$$\Psi(b)=e^{-t|D|^{2\alpha}}b^0-\int_0^te^{-(t-z)|D|^{2\alpha}}(u_a+u_b)\nabla(a+b)dz.$$
\noindent$\bullet$ Now, we introduce the space ${\bf {\mathcal X}}_T$ as follows
 $${\bf {\mathcal X}}_T=\mathcal C([0,T],{\bf{\mathcal X}}^{1-2\alpha}(\mathbb R^2))\cap L^1([0,T],{\bf{\mathcal X}}^{1}(\mathbb R^2))$$with the norm
$$\|f\|_{{\bf {\mathcal X}}_T}=\|f\|_{L^\infty_T({\bf{\mathcal X}}^{1-2\alpha})}+\|f\|_{L^1_T({\bf{\mathcal X}}^1)}.$$
Using Lemmas \ref{Lm2.4} and \ref{Lm2.5}, we can prove $\Psi({\bf {\mathcal X}}_T)\subset {\bf {\mathcal X}}_T$.\\
$\bullet$ Also, we note ${\bf B}_{r}$ the subset of ${\bf {\mathcal X}}_T$ defined by
$${\bf B}_{r}=\{\theta\in{\bf {\mathcal X}}_T;\; \|\theta\|_{L^\infty([0,T],{\bf{\mathcal X}}^{1-2\alpha})}\leq r,\;\|\theta\|_{L^1([0,T],{\bf{\mathcal X}}^{1})}\leq r\}.$$
\noindent$\bullet$ Let $b\in {\bf B}_{r}$. Prove that $\Psi(b)\in{\bf B}_{r}$. We have
$$\|\Psi(b)(t)\|_{{\bf{\mathcal X}}^{1-2\alpha}}\leq \sum_{k=0}^4I_k,$$
where $$\begin{array}{ccc}
I_0&=&\|e^{-t|D|^{2\alpha}}b^0\|_{{\bf{\mathcal X}}^{1-2\alpha}}\\
I_1&=&\displaystyle\int_0^t\|e^{-(t-z)|D|^{2\alpha}}u_a\nabla a\|_{{\bf{\mathcal X}}^{1-2\alpha}}dz\\
I_2&=&\displaystyle\int_0^t\|e^{-(t-z)|D|^{2\alpha}}u_a\nabla b\|_{{\bf{\mathcal X}}^{1-2\alpha}}dz\\
I_3&=&\displaystyle\int_0^t\|e^{-(t-z)|D|^{2\alpha}}u_b\nabla a\|_{{\bf{\mathcal X}}^{1-2\alpha}}dz\\
I_4&=&\displaystyle\int_0^t\|e^{-(t-z)|D|^{2\alpha}}u_b\nabla b\|_{{\bf{\mathcal X}}^{1-2\alpha}}dz.
\end{array}$$
Using equation (\ref{eqv01}), Lemma \ref{Lm2.4} and the fact $b\in {\bf B}_{r}$, we get
$$\begin{array}{ccc}
I_0&\leq&\displaystyle\frac{r}{5}\;\;\;\;\;\;\;\;\;\;\;\;\;\;\;\;\;\;\;\;\;\;\;\;\;\;\;\;\;\;\;\;\;\;\;\;\;\;\;\;\;\;\;\;\;\;\;\;\;\;\;\;\;\;\;\\
I_1&\leq&\displaystyle4\|a\|_{L^\infty_T({\bf{\mathcal X}}^{1-2\alpha})}\|a\|_{L^1_T({\bf{\mathcal X}}^1)}\leq 4\varepsilon\|\theta^0\|_{({\bf{\mathcal X}}^{1-2\alpha}}\leq \frac{r}{5}\\
I_2,\;I_3&\leq&2\|a\|_{L^\infty_T({\bf{\mathcal X}}^{1-2\alpha})}^{1-\frac{1}{2\alpha}}
\|a\|_{L^1_T({\bf{\mathcal X}}^1)}^{\frac{1}{2\alpha}}
\|b\|_{L^\infty_T({\bf{\mathcal X}}^{1-2\alpha})}^{\frac{1}{2\alpha}}\|b\|_{L^1_T({\bf{\mathcal X}}^1)}^{1-\frac{1}{2\alpha}}\\
&&+2\|a\|_{L^\infty_T({\bf{\mathcal X}}^{1-2\alpha})}^{\frac{1}{2\alpha}}
\|a\|_{L^1_T({\bf{\mathcal X}}^1)}^{1-\frac{1}{2\alpha}}
\|b\|_{L^\infty_T({\bf{\mathcal X}}^{1-2\alpha})}^{1-\frac{1}{2\alpha}}\|b\|_{L^1_T({\bf{\mathcal X}}^1)}^{\frac{1}{2\alpha}}\\
&\leq&2\Big(\|\theta^0\|_{\bf{\mathcal X}^{1-2\alpha}}^{1-\frac{1}{2\alpha}}
\varepsilon^{\frac{1}{2\alpha}}+\|\theta^0\|_{\bf{\mathcal X}^{1-2\alpha}}^{\frac{1}{2\alpha}}
\varepsilon^{1-\frac{1}{2\alpha}}\Big)r\\
&\leq&\displaystyle\frac{r}{5}\;\;\;\;\;\;\;\;\;\;\;\;\;\;\;\;\;\;\;\;\;\;\;\;\;\;\;\;\;\;\;\;\;\;\;\;\;\;\;\;\;\;\;\;\;\;\;\;\;\;\;\;\;\;\;\\
I_4&\leq&\displaystyle4\|b\|_{L^\infty_T({\bf{\mathcal X}}^{1-2\alpha})}\|b\|_{L^1_T({\bf{\mathcal X}}^1)}\leq 4r^2\leq\frac{r}{5}.\;\;\;\;\;\;\;\;\;\;\;\;\\
\end{array}$$
Then
\begin{equation}\label{eqv3}
\|\Psi(b)(t)\|_{{\bf{\mathcal X}}^{1-2\alpha}}\leq r.
\end{equation}

Similarly,
$$\|\Psi(b)(t)\|_{L^1_T({\bf{\mathcal X}}^{1})}\leq \sum_{k=0}^4J_k,$$
where $$\begin{array}{ccc}
J_0&=&\displaystyle\int_0^T\|e^{-t|D|^{2\alpha}}b^0\|_{{\bf{\mathcal X}}^{1}}dt,\\
J_1&=&\displaystyle\int_0^T\|\int_0^te^{-(t-z)|D|^{2\alpha}}u_a\nabla adz\|_{{\bf{\mathcal X}}^{1}}dt,\\
J_2&=&\displaystyle\int_0^T\|\int_0^te^{-(t-z)|D|^{2\alpha}}u_a\nabla bdz\|_{{\bf{\mathcal X}}^{1}}dt,\\
J_3&=&\displaystyle\int_0^T\|\int_0^te^{-(t-z)|D|^{2\alpha}}u_b\nabla adz\|_{{\bf{\mathcal X}}^{1}}dt,\\
J_4&=&\displaystyle\int_0^T\|\int_0^te^{-(t-z)|D|^{2\alpha}}u_b\nabla bdz\|_{{\bf{\mathcal X}}^{1}}dt.\\
\end{array}$$
By using Lemma \ref{Lm2.5}, we have
$$\begin{array}{ccc}
J_0&\leq&\displaystyle\frac{r}{5}\;\;\;\;\;\;\;\;\;\;\;\;\;\;\;\;\;\;\;\;\;\;\;\;\;\;\;\;\;\;\;\;\;\;\;\;\;\;\;\;\;\;\;\;\;\;\;\;\;\;\;\;\;\;\;\\
J_1&\leq&\displaystyle4\|a\|_{L^\infty_T({\bf{\mathcal X}}^{1-2\alpha})}\|a\|_{L^1_T({\bf{\mathcal X}}^1)}\leq 4\varepsilon\|\theta^0\|_{({\bf{\mathcal X}}^{1-2\alpha}}\leq \frac{r}{5}\\
J_2,\;J_3&\leq&2\|a\|_{L^\infty_T({\bf{\mathcal X}}^{1-2\alpha})}^{1-\frac{1}{2\alpha}}
\|a\|_{L^1_T({\bf{\mathcal X}}^1)}^{\frac{1}{2\alpha}}
\|b\|_{L^\infty_T({\bf{\mathcal X}}^{1-2\alpha})}^{\frac{1}{2\alpha}}\|b\|_{L^1_T({\bf{\mathcal X}}^1)}^{1-\frac{1}{2\alpha}}\\
&&+2\|a\|_{L^\infty_T({\bf{\mathcal X}}^{1-2\alpha})}^{\frac{1}{2\alpha}}
\|a\|_{L^1_T({\bf{\mathcal X}}^1)}^{1-\frac{1}{2\alpha}}
\|b\|_{L^\infty_T({\bf{\mathcal X}}^{1-2\alpha})}^{1-\frac{1}{2\alpha}}\|b\|_{L^1_T({\bf{\mathcal X}}^1)}^{\frac{1}{2\alpha}}\\
&\leq&2\Big(\|\theta^0\|_{\bf{\mathcal X}^{1-2\alpha}}^{1-\frac{1}{2\alpha}}
\varepsilon^{\frac{1}{2\alpha}}+\|\theta^0\|_{\bf{\mathcal X}^{1-2\alpha}}^{\frac{1}{2\alpha}}
\varepsilon^{1-\frac{1}{2\alpha}}\Big)r\\
&\leq&\displaystyle\frac{r}{5}\;\;\;\;\;\;\;\;\;\;\;\;\;\;\;\;\;\;\;\;\;\;\;\;\;\;\;\;\;\;\;\;\;\;\;\;\;\;\;\;\;\;\;\;\;\;\;\;\;\;\;\;\;\;\;\\
J_4&\leq&\displaystyle4\|b\|_{L^\infty_T({\bf{\mathcal X}}^{1-2\alpha})}\|b\|_{L^1_T({\bf{\mathcal X}}^1)}\leq 4r^2\leq\frac{r}{5}.\;\;\;\;\;\;\;\;\;\;\;\;\\
\end{array}$$
Then
\begin{equation}\label{eqv4}
\|\Psi(b)(t)\|_{L^1_T({\bf{\mathcal X}}^1)}\leq r.
\end{equation}
Combining equation (\ref{eqv3}) and (\ref{eqv4}), we get $\Psi(b)\in {\bf B}_{r}$ and we can deduce
\begin{equation}\label{ex-p1}
\Psi({\bf B}_{r})\subset{\bf B}_{r}.
\end{equation}
$\bullet$ Prove that  $$\|\Psi(b_1)-\Psi(b_2)\|_{{\bf {\mathcal X}}_T}\leq \frac{1}{2}\|b_1-b_2\|_{{\bf {\mathcal X}}_T},\;b_1,b_2\in{\bf B}_{r}.$$
Put $B_1$ and $B_2$ defined by
$$B_i=(\partial_{x_2}\Lambda^{-1}b_i,\; -\partial_{x_1}\Lambda^{-1}b_i).$$
We have

$$\begin{array}{ccc}
\Psi(b_1)-\Psi(b_2)&=&-\displaystyle\int_0^te^{-(t-z)|D|^{2\alpha}}\Big(B_2\nabla b_2-B_1\nabla b_1\Big)\\
&=&-\displaystyle\int_0^te^{-(t-z)|D|^{2\alpha}}\Big(B_2\nabla( b_2-b_1)+(B_2-B_1)\nabla b_1\Big)\\
\end{array}$$
and
$$\|\Psi(b_1)-\Psi(b_2)\|_{{\bf{\mathcal X}}^{1-2\alpha}}\leq K_1+K_2,$$
with
$$\begin{array}{ccc}
K_1&=&\displaystyle\|\int_0^te^{-(t-z)|D|^{2\alpha}}B_2\nabla( b_2-b_1)dz\|_{{\bf{\mathcal X}}^{1-2\alpha}},\\
K_2&=&\displaystyle\|\int_0^te^{-(t-z)|D|^{2\alpha}}(B_2-B_1)\nabla b_1dz\|_{{\bf{\mathcal X}}^{1-2\alpha}}.\\
\end{array}$$
Using Lemmas \ref{Lm2.2} and \ref{Lm2.4}, we can deduce
$$\begin{array}{ccc}
K_1&\leq&\displaystyle2\|B_2\|_{L^\infty_{T}
({\bf{\mathcal X}}^{1-2\alpha})}^{\frac{1}{2\alpha}}\|B_2\|_{L^1_T({\bf{\mathcal X}}^{1})}^{1-\frac{1}{2\alpha}}
\|b_2-b_1\|_{L^\infty_{T}({\bf{\mathcal X}}^{1-2\alpha})}^{1-\frac{1}{2\alpha}}\|b_2-b_1\|_{L^1_T({\bf{\mathcal X}}^{1})}^{\frac{1}{2\alpha}}\\
&&+\;\;2\|B_2\|_{L^\infty_{T}({\bf{\mathcal X}}^{1-2\alpha})}^{1-\frac{1}{2\alpha}}\|B_2\|_{L^1_T({\bf{\mathcal X}}^{1})}^{\frac{1}{2\alpha}}
\|b_2-b_1\|_{L^\infty_{T}({\bf{\mathcal X}}^{1-2\alpha})}^{\frac{1}{2\alpha}}\|b_2-b_1\|_{L^1_T({\bf{\mathcal X}}^{1})}^{1-\frac{1}{2\alpha}}\\
&\leq&\displaystyle2\|b_2\|_{L^\infty_{T}
({\bf{\mathcal X}}^{1-2\alpha})}^{\frac{1}{2\alpha}}\|b_2\|_{L^1_T({\bf{\mathcal X}}^{1})}^{1-\frac{1}{2\alpha}}
\|b_2-b_1\|_{L^\infty_{T}({\bf{\mathcal X}}^{1-2\alpha})}^{1-\frac{1}{2\alpha}}\|b_2-b_1\|_{L^1_T({\bf{\mathcal X}}^{1})}^{\frac{1}{2\alpha}}\\
&&+\;\;2\|b_2\|_{L^\infty_{T}({\bf{\mathcal X}}^{1-2\alpha})}^{1-\frac{1}{2\alpha}}\|b_2\|_{L^1_T({\bf{\mathcal X}}^{1})}^{\frac{1}{2\alpha}}
\|b_2-b_1\|_{L^\infty_{T}({\bf{\mathcal X}}^{1-2\alpha})}^{\frac{1}{2\alpha}}\|b_2-b_1\|_{L^1_T({\bf{\mathcal X}}^{1})}^{1-\frac{1}{2\alpha}}.
\end{array}$$ Then
$$\begin{array}{ccc}
K_1
&\leq& 2r\|b_2-b_1\|_{L^\infty_{T}({\bf{\mathcal X}}^{1-2\alpha})}^{1-\frac{1}{2\alpha}}\|b_2-b_1\|_{L^1_T({\bf{\mathcal X}}^{1})}^{\frac{1}{2\alpha}}\\
&&+\;\;2r\|b_2-b_1\|_{L^\infty_{T}({\bf{\mathcal X}}^{1-2\alpha})}^{\frac{1}{2\alpha}}\|b_2-b_1\|_{L^1_T({\bf{\mathcal X}}^{1})}^{1-\frac{1}{2\alpha}}\\
&\leq& 4r\|b_2-b_1\|_{_{{\bf {\mathcal X}}_T}}.\\
\end{array}$$
Similarly, we get
$$K_2\leq 4r\|b_2-b_1\|_{_{{\bf {\mathcal X}}_T}}.$$
Using the fact $r<1/20$, we obtain
\begin{equation}\label{ex-p2}
\|\Psi(b_2)-\Psi(b_1)\|_{{\bf {\mathcal X}}_T}\leq \frac{1}{2}\|b_2-b_1\|_{{\bf {\mathcal X}}_T},\;\;\forall\, b_1,\;b_2\in{\bf B}_{r}.
\end{equation}
Then, combining equations (\ref{ex-p1})-(\ref{ex-p2}) and the fixed point theorem, there is a unique $b\in {\bf B}_{r}$ such that $\theta=a+b$ is solution of $({\bf S}_{\alpha})$ with $\theta\in {\bf {\mathcal X}}_T(\mathbb R^2)$.\\
\subsection{Uniqueness} \label{S3.2}Let $\theta_1,\theta_2\in \mathcal C([0,T],{\bf{\mathcal X}}^{1-2\alpha}(\mathbb R^2))$ be two solution of $({\bf S}_{\alpha})$, with $\theta_1\in L^1([0,T],{\bf{\mathcal X}}^{1-2\alpha}(\mathbb R^2))$ and $\theta_1(0)=\theta_2(0)$. Put $u_1$, $u_2$, $\delta$ and $w$ defined as follows
  $$u_i=(\partial_{x_2}\Lambda^{-1}\theta_i,\; -\partial_{x_1}\Lambda^{-1}\theta_i).$$
  $$\delta=\theta_1-\theta_2$$
  $$w=u_1-u_2$$
We have
\begin{equation}\label{exiseq10}\partial_t\delta+(-\Delta)^\alpha\delta+w\nabla\delta+w\nabla \theta_1+u_1\nabla\delta=0\end{equation}
then
$$\partial_t\widehat{\delta}+|\xi|^{2\alpha}\widehat{\delta}+\mathcal F(w\nabla \delta)+\mathcal F(w\nabla \theta_1)+\mathcal F(u_1\nabla\delta)=0$$
multiply this equation by $\overline{{\widehat{\delta}}}$, we get
\begin{equation}\label{exiseq11}\partial_t\widehat{\delta}.\overline{{\widehat{\delta}}}+|\xi|^{2\alpha}|\widehat{\delta}|^2+\mathcal F(w\nabla \delta).\overline{{\widehat{\delta}}}+\mathcal F(w\nabla \theta_1).\overline{{\widehat{\delta}}}+\mathcal F(u_1\nabla\delta).\overline{{\widehat{\delta}}}=0.\end{equation}
From equation (\ref{exiseq10}), we have
  $$\partial_t\overline{\widehat{\delta}}+|\xi|^{2\alpha}\overline{\widehat{\delta}}+\overline{\mathcal F(w\nabla \delta)}+\overline{\mathcal F(w\nabla \theta_1)}+\overline{\mathcal F(u_1\nabla\delta)}=0$$
    multiply this equation by ${\widehat{\delta}}$, we get
   \begin{equation}\label{exiseq12}\partial_t\overline{\widehat{\delta}}.{{\widehat{\delta}}}+|\xi|^{2\alpha}|\widehat{\delta}|^2+\overline{\mathcal F(w\nabla \delta)}.{\widehat{\delta}}+\overline{\mathcal F(w\nabla \theta_1)}.{\widehat{\delta}}+\overline{\mathcal F(u_1\nabla\delta)}.{\widehat{\delta}}=0.\end{equation}
By summing (\ref{exiseq11}) and (\ref{exiseq12}), we get
$$\partial_t|{\widehat{\delta}}|^2+2|\xi|^{2\alpha}|\widehat{\delta}|^2+2Re\Big({\mathcal F(w\nabla \delta)}.\overline{{\widehat{\delta}}}\Big)+2Re\Big({\mathcal F(w\nabla \theta_1)}.\overline{\widehat{\delta}}\Big)+2Re\Big({\mathcal F(u_1\nabla\delta)}.\overline{{\widehat{\delta}}}\Big)=0$$
 and
  $$\partial_t|{\widehat{\delta}}|^2+2|\xi|^{2\alpha}|\widehat{\delta}|^2\leq 2|{\mathcal F(w\nabla \delta)}|.|\overline{{\widehat{\delta}}}|+2|{\mathcal F(w\nabla \theta_1)}|.|\overline{\widehat{\delta}}|+2|{\mathcal F(u_1\nabla\delta)}|.|\overline{{\widehat{\delta}}}|.$$
For $\varepsilon>0$, we have
$$\partial_t|{\widehat{\delta}}|^2=\partial_t(|{\widehat{\delta}}|^2+\varepsilon^2)
=2\sqrt{|{\widehat{\delta}}|^2+\varepsilon^2}\partial_t\sqrt{|{\widehat{\delta}}|^2+\varepsilon^2}$$
then
$$2\partial_t\sqrt{|{\widehat{\delta}}|^2+\varepsilon^2}
+2|\xi|^{2\alpha}\frac{|\widehat{\delta}|^2}{\sqrt{|{\widehat{\delta}}|^2+\varepsilon^2}}
\leq 2|{\mathcal F(w\nabla \delta)}|\frac{|{\widehat{\delta}}|}{\sqrt{|{\widehat{\delta}}|^2+\varepsilon^2}}+2|{\mathcal F(w\nabla \theta_1)}|\frac{|{\widehat{\delta}}|}{\sqrt{|{\widehat{\delta}}|^2+\varepsilon^2}}$$$$\quad\quad\quad\quad\quad\quad\quad\quad\quad\quad\quad\quad+\quad2|{\mathcal F(u_1\nabla\delta)}|\frac{|{\widehat{\delta}}|}{\sqrt{|{\widehat{\delta}}|^2+\varepsilon^2}}$$
$$\quad\quad\quad\quad\quad\quad\quad\quad\quad\quad\quad\quad\leq 2|{\mathcal F(w\nabla \delta)}|+2|{\mathcal F(w\nabla \theta_1)}|+2|{\mathcal F(u_1\nabla\delta)}|.$$
By integrating with respect to time
$$\sqrt{|{\widehat{\delta}}|^2+\varepsilon^2}
+\int_0^t|\xi|^{2\alpha}\frac{|\widehat{\delta}|^2}{\sqrt{|{\widehat{\delta}}|^2+\varepsilon^2}}
\leq \int_0^t|{\mathcal F(w\nabla \delta)}|+\int_0^t|{\mathcal F(w\nabla \theta_1)}|+\int_0^t|{\mathcal F(u_1\nabla\delta)}|.$$
Letting $\varepsilon\rightarrow0$, we get
$$|{\widehat{\delta}}|
+\int_0^t|\xi|^{2\alpha}|\widehat{\delta}|
\leq \int_0^t|{\mathcal F(w\nabla \delta)}|+\int_0^t|{\mathcal F(w\nabla \theta_1)}|+\int_0^t|{\mathcal F(u_1\nabla\delta)}|.$$
Multiplying by $|\xi|^{1-2\alpha}$ and integrating with respect to $\xi$, we get
$$\|\delta\|_{{\bf{\mathcal X}}^{1-2\alpha}}+\int_0^t\|\delta\|_{{\bf{\mathcal X}}^{1}}\leq \int_0^t\|w\nabla \delta\|_{{\bf{\mathcal X}}^{1-2\alpha}}+\int_0^t\|w\nabla \theta_1\|_{{\bf{\mathcal X}}^{1-2\alpha}}+\int_0^t\|u_1\nabla\delta\|_{{\bf{\mathcal X}}^{1-2\alpha}}.$$
$$\leq \int_0^t\|\delta w\|_{{\bf{\mathcal X}}^{2-2\alpha}}+\int_0^t\|\theta_1 w\|_{{\bf{\mathcal X}}^{2-2\alpha}}+\int_0^t\|\delta u_1\|_{{\bf{\mathcal X}}^{2-2\alpha}}.$$
By Lemma \ref{Lm2.2}, we get
$$\|\delta\|_{{\bf{\mathcal X}}^{1-2\alpha}}+\int_0^t\|\delta\|_{{\bf{\mathcal X}}^{1}}\leq 4\int_0^t\|\delta\|_{{\bf{\mathcal X}}^{1-2\alpha}}\|\delta\|_{{\bf{\mathcal X}}^{1}}+
2\int_0^t\|\delta\|_{{\bf{\mathcal X}}^{1-2\alpha}}^{1-\frac{1}{2\alpha}}\|\delta\|_{{\bf{\mathcal X}}^{1}}^{\frac{1}{2\alpha}}
\|\theta_1\|_{{\bf{\mathcal X}}^{1-2\alpha}}^{\frac{1}{2\alpha}}\|\theta_1\|_{{\bf{\mathcal X}}^{1}}^{1-\frac{1}{2\alpha}}
$$$$+2\int_0^t\|\delta\|_{{\bf{\mathcal X}}^{1-2\alpha}}^{\frac{1}{2\alpha}}\|\delta\|_{{\bf{\mathcal X}}^{1}}^{1-\frac{1}{2\alpha}}
\|\theta_1\|_{{\bf{\mathcal X}}^{1-2\alpha}}^{1-\frac{1}{2\alpha}}\|\theta_1\|_{{\bf{\mathcal X}}^{1}}^{\frac{1}{2\alpha}}.$$
Using the elementary inequality
$$xy\leq \frac{x^p}{p}+\frac{y^q}{q}$$
with
$$p=\frac{1}{2\alpha},\;q=\frac{2\alpha}{2\alpha-1},$$
we get
$$\|\delta\|_{{\bf{\mathcal X}}^{1-2\alpha}}^{1-\frac{1}{2\alpha}}\|\delta\|_{{\bf{\mathcal X}}^{1}}^{\frac{1}{2\alpha}}
\|\theta_1\|_{{\bf{\mathcal X}}^{1-2\alpha}}^{\frac{1}{2\alpha}}\|\theta_1\|_{{\bf{\mathcal X}}^{1}}^{1-\frac{1}{2\alpha}}\leq \frac{1}{4}\|\delta\|_{{\bf{\mathcal X}}^{1}}
+c_\alpha \|\delta\|_{{\bf{\mathcal X}}^{1-2\alpha}}\|\theta_1\|_{{\bf{\mathcal X}}^{1-2\alpha}}^{\frac{1}{2\alpha-1}}\|\theta_1\|_{{\bf{\mathcal X}}^{1}}$$
and
$$\|\delta\|_{{\bf{\mathcal X}}^{1-2\alpha}}^{\frac{1}{2\alpha}}\|\delta\|_{{\bf{\mathcal X}}^{1}}^{1-\frac{1}{2\alpha}}
\|\theta_1\|_{{\bf{\mathcal X}}^{1-2\alpha}}^{1-\frac{1}{2\alpha}}\|\theta_1\|_{{\bf{\mathcal X}}^{1}}^{\frac{1}{2\alpha}}\leq \frac{1}{4}\|\delta\|_{{\bf{\mathcal X}}^{1}}
+c'_\alpha \|\delta\|_{{\bf{\mathcal X}}^{1-2\alpha}}\|\theta_1\|_{{\bf{\mathcal X}}^{1-2\alpha}}^{\frac{1}{2\alpha-1}}\|\theta_1\|_{{\bf{\mathcal X}}^{1}}.$$
Then
$$\|\delta\|_{{\bf{\mathcal X}}^{1-2\alpha}}+\int_0^t\|\delta\|_{{\bf{\mathcal X}}^{1}}\leq (c_\alpha+c'_\alpha)\int_0^t \|\delta\|_{{\bf{\mathcal X}}^{1-2\alpha}}\|\theta_1\|_{{\bf{\mathcal X}}^{1-2\alpha}}^{\frac{1}{2\alpha-1}}\|\theta_1\|_{{\bf{\mathcal X}}^{1}}.$$
Using Gronwall Lemma and the fact $(t\mapsto\|\theta_1\|_{{\bf{\mathcal X}}^{1-2\alpha}}^{\frac{1}{2\alpha-1}}\|\theta_1\|_{{\bf{\mathcal X}}^{1}})\in L^1([0,T])$, we can deduce that $\delta=0$ in $[0,T]$ which gives the uniqueness.
\subsection{Small initial data}\label{S3.3}
In this section, we suppose that $\|\theta^0\|_{{\bf{\mathcal X}}^{1-2\alpha}}<1/4.$\\ Pass to Fourier transform of the first equation of $({\bf S}_\alpha)$ and  multiply it by $\overline{{\widehat{\theta}}}$, we get
\begin{equation}\label{smleq1}\partial_t\widehat{\theta}.\overline{{\widehat{\theta}}}+|\xi|^{2\alpha}|\widehat{\theta}|^2+\mathcal F(u\nabla \theta).\overline{{\widehat{\theta}}}=0.\end{equation}
Similarly, we have
 $$\partial_t\overline{\widehat{\theta}}+|\xi|^{2\alpha}\overline{\widehat{\theta}}+\overline{\mathcal F(u\nabla \theta)}=0$$
and
  \begin{equation}\label{smleq2}\partial_t\overline{\widehat{\theta}}.{{\widehat{\theta}}}+|\xi|^{2\alpha}|\widehat{\theta}|^2+\overline{\mathcal F(u\nabla \theta)}.{\widehat{\theta}}=0.\end{equation}
By summing (\ref{smleq1}) and (\ref{smleq2}), we get
$$\partial_t|{\widehat{\theta}}|^2+2|\xi|^{2\alpha}|\widehat{\theta}|^2+2Re\Big({\mathcal F(u\nabla \theta)}.\overline{{\widehat{\theta}}}\Big)=0$$
 and
  $$\partial_t|{\widehat{\theta}}|^2+2|\xi|^{2\alpha}|\widehat{\theta}|^2\leq 2|{\mathcal F(u\nabla \theta)}|.|\overline{{\widehat{\theta}}}|.$$
For $\varepsilon>0$, we have
$$\partial_t|{\widehat{\theta}}|^2=\partial_t(|{\widehat{\theta}}|^2+\varepsilon^2)=2\sqrt{|{\widehat{\theta}}|^2+\varepsilon^2}
\partial_t\sqrt{|{\widehat{\theta}}|^2+\varepsilon^2}.$$
Then
$$2\partial_t\sqrt{|{\widehat{\theta}}|^2+\varepsilon^2}
+2|\xi|^{2\alpha}\frac{|\widehat{\theta}|^2}{\sqrt{|{\widehat{\theta}}|^2+\varepsilon^2}}
\leq 2|{\mathcal F(u\nabla \theta)}|.\frac{|{\widehat{\theta}}|}{\sqrt{|{\widehat{\theta}}|^2+\varepsilon^2}}$$
$$\quad\quad\quad\quad\quad\quad\quad\quad\quad\quad\leq 2|{\mathcal F(u\nabla \theta)}|.$$
By integrating with respect to time
$$\sqrt{|{\widehat{\theta}}|^2+\varepsilon^2}
+\int_0^t|\xi|^{2\alpha}\frac{|\widehat{\theta}|^2}{\sqrt{|{\widehat{\delta}}|^2+\varepsilon^2}}
\leq \sqrt{|{\widehat{\theta^0}}|^2+\varepsilon^2}+\int_0^t|{\mathcal F(w\nabla \delta)}|.$$
Letting $\varepsilon\rightarrow0$, we get
$$|{\widehat{\theta}}|
+\int_0^t|\xi|^{2\alpha}|\widehat{\theta}|
\leq |{\widehat{\theta^0}}|+\int_0^t|{\mathcal F(u\nabla \theta)}|.$$
Multiplying by $|\xi|^{1-2\alpha}$ and integrating with respect to $\xi$, we get
$$\|\theta\|_{{\bf{\mathcal X}}^{1-2\alpha}}+\int_0^t\|\theta\|_{{\bf{\mathcal X}}^{1}}\leq \|\theta^0\|_{{\bf{\mathcal X}}^{1-2\alpha}}+\int_0^t\|u\nabla \theta\|_{{\bf{\mathcal X}}^{1-2\alpha}}.$$
By lemma \ref{Lm2.2}, we get
$$\|\theta\|_{{\bf{\mathcal X}}^{1-2\alpha}}+\int_0^t\|\theta\|_{{\bf{\mathcal X}}^{1}}\leq \|\theta^0\|_{{\bf{\mathcal X}}^{1-2\alpha}}+4\int_0^t\|\theta\|_{{\bf{\mathcal X}}^{1-2\alpha}}\|\theta\|_{{\bf{\mathcal X}}^{1}}$$
$$\leq \|\theta^0\|_{{\bf{\mathcal X}}^{1-2\alpha}}+4\sup_{z\in[0,t]}\|\theta(z)\|_{{\bf{\mathcal X}}^{1-2\alpha}}\int_0^t\|\theta\|_{{\bf{\mathcal X}}^{1}}.$$
Let $T=\sup\{t>0;\;\sup_{z\in[0,t]}\|\theta(z)\|_{{\bf{\mathcal X}}^{1-2\alpha}}<\frac{1}{4}\}.$ By the above equation, we have
$$\|\theta(t)\|_{{\bf{\mathcal X}}^{1-2\alpha}}\leq \|\theta^0\|_{{\bf{\mathcal X}}^{1-2\alpha}}< \frac{1}{4},\forall \, t\in[0,T),$$
then $T=\infty$. Therefore the global existence and inequality (\ref{thjm1eq}) are proved.
\section {Proof of Theorem \ref{thjm2} }\label{S4}The same approach of Theorem \ref{thjm1} is used to obtain a blow up result of $\theta\in\mathcal C([0,T^*),{\bf{\mathcal X}}^{1-2\alpha}(\mathbb R^2))$ if $T^*<\infty$.\\
Assume that $\displaystyle\int_0^{T^*}\|\theta(t)\|_{{\bf{\mathcal X}}^{1}} dt<\infty.$ Let $0<T_0<T^*$ such that
$$\int_{T_0}^{T^*}\|\theta(t)\|_{{\bf{\mathcal X}}^{1}} dt<\frac{1}{2}.$$
For $t\in[T_0,T^*)$ and $s\in[T_0,t]$
$$\|\theta(s)\|_{{\bf{\mathcal X}}^{1-2\alpha}}+\int_{T_0}^s\|\theta(t)\|_{{\bf{\mathcal X}}^{1}}\leq
\|\theta(T_0)\|_{{\bf{\mathcal X}}^{1-2\alpha}}+\int_{T_0}^s\|\theta(t)\|_{{\bf{\mathcal X}}^{1-2\alpha}}\|\theta(t)\|_{{\bf{\mathcal X}}^{1}}$$
$$\leq
\|\theta(T_0)\|_{{\bf{\mathcal X}}^{1-2\alpha}}+\sup_{T_0\leq z\leq t}\|\theta(z)\|_{{\bf{\mathcal X}}^{1-2\alpha}}\int_{T_0}^t\|\theta(t)\|_{{\bf{\mathcal X}}^{1}}$$
$$\leq
\|\theta(T_0)\|_{{\bf{\mathcal X}}^{1-2\alpha}}+\frac{1}{2}\sup_{T_0\leq z\leq t}\|\theta(z)\|_{{\bf{\mathcal X}}^{1-2\alpha}}$$
then
$$\sup_{0\leq s\leq t}\|\theta(s)\|_{{\bf{\mathcal X}}^{1-2\alpha}}\leq
\|\theta(T_0)\|_{{\bf{\mathcal X}}^{1-2\alpha}}+\frac{1}{2}\sup_{T_0\leq z\leq t}\|\theta(z)\|_{{\bf{\mathcal X}}^{1-2\alpha}}.$$
We can deduce that
 $$\sup_{0\leq s\leq t}\|\theta(s)\|_{{\bf{\mathcal X}}^{1-2\alpha}}\leq2
\|\theta(T_0)\|_{{\bf{\mathcal X}}^{1-2\alpha}},\;\;\forall t\in[T_0,T^*).$$
Put $$M=\max\Big(2\|\theta(T_0)\|_{{\bf{\mathcal X}}^{1-2\alpha}};\;\;\max_{t\in[0,T_0]}\|\theta(t)\|_{{\bf{\mathcal X}}^{1-2\alpha}}\Big).$$
We have
$$\|\theta(t)\|_{{\bf{\mathcal X}}^{1-2\alpha}}\leq M,\;\;\forall t\in[0,T^*).$$
Using the integral form of $\theta$, we can write, for $t<t'\in[0,T^*)$,
$$\theta(t')-\theta(t)=L_1(t,t')+L_2(t,t')$$
with
$$\begin{array}{ccc}
L_1&=&\displaystyle\int_0^t\Big(1-e^{-(t'-t)|D|^{2\alpha}}\Big)e^{-(t-z)|D|^{2\alpha}}(u.\nabla \theta)(z)\;dz,\\
L_2&=&\displaystyle\int_t^{t'}e^{-(t'-z)|D|^{2\alpha}}(u.\nabla \theta)(z)\;dz.\\
\end{array}$$
We have
$$\begin{array}{ccc}
\|L_1(t,t')\|_{{\bf{\mathcal X}}^{1-2\alpha}}&\leq&\displaystyle\int_0^t\int_\xi\Big(1-e^{-(t'-t)|\xi|^{2\alpha}}\Big)e^{-(t-z)|\xi|^{2\alpha}}|\xi|^{1-2\alpha}|\mathcal F(u.\nabla \theta)(z,\xi)|\;d\xi\; dz\\
&\leq&\displaystyle\int_0^t\int_\xi\Big(1-e^{-(t'-t)|\xi|^{2\alpha}}\Big)|\xi|^{1-2\alpha}|\mathcal F(u.\nabla \theta)(z,\xi)|\;d\xi\;dz.\\
\|L_2(t,t')\|_{{\bf{\mathcal X}}^{1-2\alpha}}&\leq&\displaystyle\int_t^{t'}\int_\xi e^{-(t'-z)|\xi|^{2\alpha}}|\xi|^{1-2\alpha}|\mathcal F(u.\nabla \theta)(z,\xi)|\;d\xi\;dz\\
&\leq&\displaystyle\int_t^{t'}\int_\xi |\xi|^{1-2\alpha}|\mathcal F(u.\nabla \theta)(z,\xi)|\;d\xi\;dz.
\end{array}$$
Using Lemmas \ref{Lm2.2}-\ref{Lm2.4} and Dominated Convergence Theorem, we can deduce that
$$\limsup_{\begin{array}{c}
t,t'\nearrow  T^*\\
t<t'
\end{array}}\|L_1(t,t')\|_{{\bf{\mathcal X}}^{1-2\alpha}}=0,\;\;\limsup_{\begin{array}{c}
t,t'\nearrow  T^*\\
t<t'
\end{array}}\|L_2(t,t')\|_{{\bf{\mathcal X}}^{1-2\alpha}}=0.$$
Therefore
$$\limsup_{\begin{array}{c}
t,t'\nearrow  T^*\\
t<t'
\end{array}}\|\theta(t')-\theta(t)\|_{{\bf{\mathcal X}}^{1-2\alpha}}=0.$$
Then $\theta(t)$ is of Cauchy type at the left of $T^*$ in the Banach space ${\bf{\mathcal X}}^{1-2\alpha}(\mathbb R^2)$. Then, there is $\theta^*$ an element of ${\bf{\mathcal X}}^{1-2\alpha}(\mathbb R^2)$ such that
$$\lim_{t\nearrow T^*}\theta(t)=\theta^*.$$
Now, consider the following system
$$\left\{\begin{array}{l}
  \displaystyle\partial_t
a +(-\Delta)^\alpha a+u_a.\nabla a =0 \\
a(0)=\theta^*
\end{array}\right.$$
By Theorem \ref{thjm1}, there a time $t_0>0$ and unique solution $a$ such that $a\in\mathcal C([0.t_0],{\bf{\mathcal X}}^{1-2\alpha}(\mathbb R^2))$.
Then
$$\Theta(t)=\left\{\begin{array}{l}
\theta(t),\;{\rm if}\;t\in[0,T^*)\\
a(t-T^*),\;{\rm if}\;t\in[T^*,T^*+t_0]\\
\end{array}\right.$$
is a solution of $({\bf S}_\alpha)$ with initial data $\theta^0$ on the interval $[0,T^*+t_0]$  which contradicts the maximality of $T^*$.

\end{document}